\documentclass[11pt]{amsart}
\usepackage{amsmath,graphics}

\theoremstyle{plain}
\newtheorem{theorem}{Theorem}[section]
\newtheorem{lemma}[theorem]{Lemma}

\newtheorem{proposition}[theorem]{Proposition}
\newtheorem{ther}{Theorem}

\theoremstyle{remark}
\newtheorem*{remark}{Remark}

\theoremstyle{definition}

\def \R {\mathbf{R}}
\def \Z {\mathbf{Z}}
\def \C {\mathbf{C}}
\def \H {\mathbf{H}}
\def \A {\mathbf{A}}
\def \L {\mathbf{L}}
\def \G {\mathbf{G}}
\def \F {\mathbf{F}}

\def \D {\boldsymbol{D}}
\def\bGamma{\boldsymbol{\Gamma}}

\def\Qbar{\overline{Q}}
\def\overSW{\overline{\overline{SW}}}

\def\AA{\mathcal{A}}

\def\GG{\mathcal{G}}
\def\II{\mathcal{I}}
\def\UU{\mathcal{U}}
\def\CC{\mathcal{C}}
\def\RR{\mathcal{R}}
\def\KK{\mathcal{K}}
\def\MM{\mathcal{M}}

\def\Tt{\widetilde{T}}

\def\zz{\ifmmode{z}\else{$z$}\fi}

\DeclareMathOperator{\SU}{SU}

\def\su{\ifmmode{\SU(2)}\else{$\SU(2)$}\fi}

\DeclareMathOperator{\U}{U}

\DeclareMathOperator{\spinc}{Spin^c}
\def\sw{Seiberg-Witten}
\def\SW{\ifmmode{\text{SW}}\else{$\text{SW}$}\fi}
\def\SWW{\ifmmode{\underline{\text{SW}}}\else{$\underline{\text{SW}}$}\fi}
\def\swi{Seiberg-Witten invariant}
\def\swe{Seiberg-Witten equations}

\def\cal{\mathcal}

\begin{document}

\baselineskip.525cm
\title{Mod 2 Seiberg-Witten invariants of homology tori}
\thanks{The first author was partially
supported by NSF Grant 4-50645.}
\author[Daniel Ruberman]{Daniel Ruberman}
\author[Sa\v so Strle]{Sa\v so Strle}
\address{Department of Mathematics\newline\indent
Brandeis University \newline\indent
Waltham, MA 02254}
\date\today
\maketitle

\section{Introduction}
While the \sw\ equations are defined for any $\spinc$ structure on a smooth
$4$-manifold, there is particular interest in the \sw\ equations
associated to
a spin structure.  One reason for this is that the $4$-dimensional
spin representation has a
quaternionic structure, which gives rise to a large symmetry group of
the `trivial'
reducible solution to the \sw\ equations.  This group, denoted
$J$, is generated by $\U(1)$ and the quaternion
$j$.  The interplay of this symmetry group and the deformation theory of the
trivial solution is central to Furuta's
proof of the `$10/8$' inequality, which constrains the homotopy type
of smooth spin
manifolds.  A closely related argument (known to  Kronheimer and
Furuta as well) was used by
Morgan-Szab\'o~\cite{morgan-szabo:mod2} to determine the mod
$2$ \sw\ invariant of homotopy K3 surfaces and other simply-connected
spin manifolds.

In this paper we show that the mod $2$ \sw\ invariant can be determined for a
spin manifold $X$ which has the same homology {\em groups} as the $4$-torus
$T^4$.  The value depends on the structure of the cohomology {\em ring} of $X$,
and in particular on the
$4$-fold cup product $\Lambda^4 H^1(X) \to H^4(X)$.   For the rest of the
paper, $X$ will denote a (spin) homology torus, by which we mean an
oriented spin
$4$-manifold with $H_1(X;\Z) \cong \Z^4$ and $H_2(X;\Z) \cong\Z^6$.
The cup product on
$H^2(X)$ is readily seen to be hyperbolic, but the cup product on
$H^1(X)$ is not determined
by the dimensions of these groups.  Let us define $\det(X)$, the {\em
determinant} of
$X$, to be the
absolute value of
$$
<\alpha_1\cup \alpha_2 \cup \alpha_3 \cup \alpha_4, [X]>
$$
where $\{\alpha_j\}$ is a basis for $H^1(X;\Z)$.
\begin{ther}\label{swtorus}
The value of the \swi\ for the $\spinc$ structure on $X$ with trivial
determinant line is congruent (mod $2$) to the determinant of $X$.
\end{ther}

Let $X$ be a homology torus, and let $W \to X$ be a $\spinc$ bundle
with trivial determinant line $L \to X$.  We fix a square root $L^{1/2}$ of the
determinant bundle, or equivalently a spin structure on $X$.   A
trivialization of $L^{1/2}$ provides us with a preferred  origin in
the space
$\AA$ of
$U(1)$-connections on $L^{1/2}$, namely the  (smooth)  product connection
$A_0$.   Note that although there are many spin structures on $X$,
they are all isomorphic as $\spinc$ structures; the choice of spin structure
is reflected only in the (gauge equivalence class of the) above
mentioned trivialization and hence in the choice of $A_0$, but does not
affect the argument.

Recall that the configuration space for the \swe\  is $\CC=\AA \oplus
\Gamma(W^+)$; however,  we restrict the equations to the slice
$\CC'=\KK \oplus \Gamma(W^+)$ where $\KK=\{A \in \AA|\,d^*(A-A_0)=0\}$.
With this restriction, the  moduli space $\MM$ of solutions to the
\swe\ is the quotient of the space of  solutions by the action of
the group of harmonic gauge transformations; denote by $\GG_0$
the {\em based} gauge group of harmonic gauge transformations. The formal
dimension of the moduli space $\MM$ for the trivial $\spinc$ structure on $X$
is zero (as is the index  of   Dirac operator $D_A \colon \Gamma(W^+) \to
\Gamma(W^-)$ for  any $\spinc$ connection $\nabla_A$ on $W$).  In the  absence
of any perturbation terms, however, the moduli space is not cut out
transversally because it contains the `dual' 4-torus $T^*\cong H^1(X;S^1)$  of
reducible  solutions $[A,0]$ with $A$ harmonic. This dual torus is covered
(under the action of $\GG_0$) by the  space of harmonic connections $\Tt^*=
A_0+i{\cal H}^1(X) \subset \KK$. In the case $X=T^4$ the dual torus coincides
with the moduli space; moreover, the  only point along $T^*$ where the Dirac
operator has nontrivial kernel is
$[A_0,0]$. From this, and the structure of the  quadratic term of the
equations, one can show that the \sw\ invariant (for the trivial
$\spinc$ structure) of $T^4$ is $\pm1$; this of course implies
Theorem~\ref{swtorus}  for $X=T^4$. To prove it in the general case,
we perturb the equations along $\Tt^*$ so that they become as nondegenerate
as possible; this is done in two stages - we deal with the
linear perturbations in section 2, and with nonlinear ones in section 3.
The invariant is then determined by the count of solutions which lie near
$T^*$; as  in~\cite{morgan-szabo:mod2}, we  use the involution
$j$ on the moduli space of solutions, induced by taking the dual
$\spinc$ structure,  to pair off the solutions away from $T^*$.

The spaces of sections  are $L_2^2$ unless stated otherwise; in
particular, this holds for the configuration space $\CC'$.  The gauge
transformations are in the space $L^2_3$.

\section{Dirac operators along ${T^*}$}
Recall that the \swe\ are given by a map
$SW \colon \KK \oplus \Gamma(W^+)  \to i\Omega^2_+(X) \oplus \Gamma(W^-)$,
$(A, \psi) \mapsto (F_A^+-q(\psi),D_A\psi)$, where $F_A^+$ is the self-dual
part of the curvature and $q$ is a quadratic map.  Thus for any
$A \in \Tt^*$, the linearization of the equations at $(A,0)$ is $(d^+,D_A)$.
Since
$d^+$ does not depend on $A$,
the behavior of the linearizations of  the \swe\ along $\Tt^*$ is
described by the family of Dirac operators
$\{ D_A,\ A \in \Tt^*\}$.  We think of this family as a morphism
of trivial bundles $\Tt^* \times \Gamma(W^+)
\to \Tt^* \times \Gamma(W^-)$.  Note that the gauge group $\GG_0$ acts
freely on the base space $\Tt^*$ of these bundles, so dividing out by the action
of $\GG_0$ produces bundles ${\bGamma}(W^{\pm}) \to T^*$  with fibres
$\Gamma(W^{\pm})$. Each of these bundles supports a free $J$-action which is
compatible with the $J$-action on the base $T^*$; we call any such bundle a
{\em $J$-bundle}.  The family of Dirac operators defines a family of
Fredholm operators parametrized by $T^*$; we denote the resulting
morphism by $\D \colon \bGamma(W^{+}) \to \bGamma(W^{-})$.
   From above we know  that the  pointwise index of $\D$ is $0$, but the
index {\em bundle} $Ind(\D)$ of $\D$ may be nontrivial. We will see
that the latter is determined by the cup product on $H^1(X)$.
It is, therefore, the index computation that links the cup product
structure to the behavior of the linear part of the \sw\ equations
along $\Tt^*$.

The calculation of the family index of $\D$ is similar to the one
arising in the proof~\cite{li-liu:wall,salamon:swbook} of the wall-crossing
formula for
$4$-manifolds with $b^1 > 0$. The Chern character of the index bundle
is given by $\text{ch}(Ind(\D))=\text{ch}(\L)/[X]$, since the
$\hat{A}$-genus of  $X$ is 1. Here $\L \to X  \times T^*$ is the universal
line bundle  equipped with a connection $\A$  as follows.
Let $\alpha_1, \ldots , \alpha_4 \in {\cal H}^1(X;\Z)$ be a basis and let
$t_k \mapsto 2\pi it_k\alpha_k$ be coordinates on
$T^* \cong {\cal H}^1(X;i\R)/{\cal H}^1(X;2\pi i\Z)$. The connection
1-form  of $\A$ is given by
$$2\pi i \sum_k t_k \alpha_k \, .$$
By Chern-Weil theory, the first Chern class of $\L$ is then represented
by the 2-form
$$\Omega = \sum_k  \alpha_k \wedge dt_k$$
and therefore
$$ \text{ch}(\L)=1+\Omega+\frac{1}{2}\Omega^2+ \frac{1}{6}\Omega^3 +
\frac{1}{24}\Omega^4 \, .$$
It is at this point that the cup product structure of $X$ shows up;
the formula for the Chern character of the index bundle gives
$$ \text{ch}(Ind(\D))=\pm r  [vol_{T^*}]$$
where $r$ denotes the determinant of $X$.  Consequently, $c_2(Ind(\D))=\pm r$
and $c_1(Ind(\D))=0$; this suggests that there is a simple model for the index
bundle and the construction of this model occupies the rest of the section.

In the proposition below we construct a generic model
for the index bundle, realized  by stabilizing the domain and the range of
the operators. This corresponds to a stabilization of the \swe; we will
define the stabilized equations via a map
$$\overline{SW} \colon \KK \oplus \Gamma(W^+) \oplus \H^n \to
i\Omega^2_+(X) \oplus \Gamma(W^-) \oplus\H^n\, . $$
\begin{proposition}\label{linear} Suppose that $\det(X) =r$. Then there is a
$J \times \GG_0$--equivariant stabilization of the \sw\ equations, with
reducible solutions along $T^*$, such that the  corresponding family of Dirac
operators has nontrivial kernel at exactly $r$  points on $T^*$.
\end{proposition}
\begin{proof}
As an element of the ordinary $K$-theory of $T^*$, the index bundle may be
represented as a difference of complex vector bundles.  We need to represent
$Ind(\D)$ as the difference of two genuine $J$-bundles over $T^*$, and so adapt
the standard argument to the context of $J$-equivariant $K$-theory
(compare~\cite{atiyah:reality,iriye:k-theory}) as follows.  By standard
arguments there exists a $\C$-linear  morphism $\G_0 \colon T^*
\times \C^n \to
\bGamma(W^-)$  which is onto the cokernel of $\D$. This morphism
extends to a $J$-equivariant morphism $\G \colon T^* \times \H^n  \to
\bGamma(W^-)$, where the product bundle $T^* \times \H^n$ has
the product action of $j$; the $j$-action on the space of
quaternions  $\H=\C \oplus j\C$ is  via right quaternionic
multiplication. Given any $[A,0] \in T^*$ and $w \in \C^n$ set
$$\G([A,0],jw):=j\cdot \G_0([j(A),0], \overline w)$$
and extend by linearity.  Perturbing the family of Dirac operators $\D$
by the morphism $\G$ produces a $J$-equivariant epimorphism
\begin{equation*}
\begin{align*}
\overline\D  \colon \bGamma(W^+) \oplus T^* \times \H^n  &\to
\bGamma(W^-)\\ ([A,0],\psi,w)  &\mapsto  \D_{[A,0]}\psi\,+ \G([A,0],w)\, .
\end{align*}
\end{equation*}
Through this we have represented $Ind(\D)$ in $J$-equivariant K-theory as the
difference of the kernel bundle of $\overline\D$ and the product $J$-bundle
$T^* \times \H^n$.

Considered as a complex bundle, the kernel bundle  $K:= \ker
\overline\D$ splits as a sum $K=K' \oplus K''$, for dimensional reasons, where
$K'$ is a trivial  complex bundle, and $K''$ is a $\C^2$-bundle with
$c_2(K'')=\pm r$ and $c_1(K'')=0$.   In fact, $K$ splits in the category of
$J$-bundles over $T^*$, in such a way that $K'$ is a trivial $\H^{n-1}$-bundle.
To construct this splitting note that any
vector bundle over $T^*$ with fibre dimension greater than 4 (over $\R$)
admits a nowhere vanishing section $s$. On any $J$-bundle $M \to T^*$
such a section $s$ gives rise to a trivial $J$-invariant  sub-bundle
$N \to T^*$ of complex  rank 2, spanned by $s$ and $\overline s([A,0])=
j\cdot s([j(A),0])$. Moreover, $N$ has a $J$-invariant complement
in $M$; the latter can be taken to be perpendicular to $N$ with respect
to some (compatible) hermitian inner product on $M$. For the case at
hand we choose the standard
hermitian structure on $T^* \times \H^n$ and the $L^2$-inner  product on the
fibres of $\bGamma(W^+)$. Denote the resulting $J$-equivariant isomorphism by
$\F' \colon   K' \to T^* \times \H^{n-1}$.

The bundle $K''$ admits a structure of a quaternionic line bundle; we use
this to construct a $J$-equivariant morphism $\F'' \colon K'' \to
T^* \times \H$, injective everywhere except at $r$ chosen  points on
$T^*$. Let $\RR \subset T^*$ be a $j$-invariant subset with $r$ elements;
such exists for any $r$ since $j$-action on $T^*$ has  fixed points
(for  example $[A_0,0]$). We choose a section $s_0$ of the bundle $K''$
which vanishes only at the points of $\RR$ and intersects the
zero section transversely. Then the sections $s_0$ and $\overline s_0$
(defined from $s_0$ as above) endow $K''$ with a structure of a quaternionic
line bundle over the complement of $\RR$. Dividing $s_0$ by the square of
its (quaternionic) norm produces a nowhere vanishing section $s$ of $K''$
over the complement of $\RR$ and this section $s$ induces the required
bundle morphism $\F''$.

Note that close to any $[A_k,0] \in \RR$, the norms of
linear maps $\F_{[A,0]}^{\prime\prime}$ are bounded below by some positive
constant times distance from $[A,0]$ to $[A_k,0]$. The morphisms $\F'$ and
$\F''$ together define a $J$-equivariant morphism $\F \colon K \to T^*
\times \H^n$ which is injective on all the fibres except over the points
of  $\RR$ where the kernels can be identified with a copy of $\H$.
We think of the pair $(\overline\D,\F)$ as the family of Dirac operators
associated to the stabilized \swe\ (which are defined below). Note that
by construction of $\overline\D$ and $\F$, the associated family of Dirac
operators has nontrivial kernels only at the points of $\RR$, thus proving
the last statement of the proposition.

To finish the construction of the stabilized equations, we need to globalize
the perturbation terms $\G$ and $\F$. Let $P \colon \KK \to \Tt^*$ be the
$L^2$-orthogonal projection (where we treat
$A_0$ as the  origin of the above affine spaces), $Q \colon \KK \to
(\Tt^*)^\perp$ the orthogonal  projection to the complement, and
$\Pi \colon \bGamma(W^+) \oplus T^*  \times \H^n \to  K$ the orthogonal
projection to the kernel of $\overline \D$. The morphism $\F$ defines a map
$$ F \colon   \KK \oplus \Gamma(W^+) \oplus \H^n \to \H^n$$
given by $F(A,\psi,w)= pr_2 \circ \F(\Pi([P(A),\psi,w])).$
Similarly, $\G$ gives rise to
$$ G \colon \KK \oplus \Gamma(W^+) \oplus \H^n \to \Gamma(W^-) $$
which is well defined up to gauge change by $[P(A),G(A,\psi,w)]=
\G([P(A),0],w)$; it is completely determined by the appropriate choice
of $G_{A_0}$.

We define the stabilized \swe\ via a map
$$\KK \oplus \Gamma(W^+) \oplus \H^n \to i\Omega^2_+(X)
\oplus \Gamma(W^-) \oplus \H^n $$
which is the sum of the original \sw\ map and the stabilization term
given by
$$ (A,\psi,w) \mapsto  \beta(Q(A),\psi,w)\cdot
(0,G(A,\psi,w),F(A,\psi,w))+(1-\beta(Q(A),\psi,w)) \cdot (0,0,w)$$
where $\beta$ depends smoothly on the $L_2^2$-norms of $A$ and
$\psi$ and on the norm of $w$ in such a way that it is equal to 1
in a small neighborhood of $(0,0,0)$ and equal to 0 in a slightly
bigger neighborhood; notice that $\beta(Q(-),-,-)$ is invariant
under the action of the gauge group $\GG_0$ as well as under the
action of $J$. It is clear from the nature of the perturbation terms
that the moduli space of solutions to the stabilized \sw\ equation
$\overline{SW}=0$ still contains  the torus of reducibles $T^*$.
This proves the proposition.
\end{proof}
\begin{remark}
The proof of the proposition implies not only that $T^*$ is contained
in the moduli space of solutions to $\overline{SW}=0$, but also that
it is isolated, at least away from the points of $\RR$. More precisely,
for any neighborhood $U$ of $\RR$, the complement $T^* \setminus U$ is
isolated in the moduli space.  This follows from the fact that $\F$ is
injective on the kernels of the perturbed family of Dirac operators
$\overline \D$ (along $T^*$) away from $\RR$.
\end{remark}

\section{Kuranishi maps at the points of $\RR$}
In this section we will construct a further perturbation of the stabilized
\sw\ map $\overline{SW}$ whose solution space has a particularly simple
form in a neighborhood of $T^*$, as described in the proposition below.
The perturbation is supported in a  neighborhood of the set $\RR \subset
T^*$ and is constructed by modifying the Kuranishi maps at the points
of $\RR$; these are the only points on $T^*$ at which the stabilized Dirac
operators have nontrivial kernels.
\begin{proposition}
There exists a $J \times \GG_0$-equivariant perturbation of the stabilized
\swe, such that the perturbed map $\overSW$ satisfies the following:
\begin{enumerate}
\item The torus of reducibles $T^*$ is contained and isolated in the moduli
space of solutions to the perturbed equations $\overSW=0$.
\item Given a small generic $\omega \in i\Omega^2_+(X)$ there exists an
invariant neighborhood $\UU$ of $\Tt^*$, such that all the solutions to
$\overSW=(\omega,0,0)$ that lie in $\UU$ are smooth and irreducible.
More precisely, every point in $\RR$ gives rise to a smooth circle of
solutions to $\overSW=(\omega,0,0)$ in $\UU$, contributing $\pm 1$ to the
invariant, and there are no other solutions in $\UU$.
\end{enumerate}
\end{proposition}
\begin{proof}
Points of $\RR$  fall into two categories depending on whether they are
$j$-fixed or not. We consider the former case first, making use of the
$j$-equivariance of the Kuranishi map. Then we modify the argument to deal
with the rest of the points in $\RR$.

Suppose $[A_k,0] \in \RR$ is $j$-fixed. The Kuranishi model for the
solutions  to $\overline{SW}=0$ around $(A_k,0,0)$ is given by a
$J$-equivariant map $Q \colon \R^4 \oplus \H \to \R^3 \oplus \H$, where
$\R^4$ corresponds to the harmonic 1-forms, $\R^3$ to the  self-dual
harmonic 2-forms, and the quaternions represent the kernel and the
cokernel of the perturbed Dirac operator at $A_k$. Note that the leading
term of $Q$ is a quadratic polynomial map which we will make non-degenerate
by a perturbation. Denote by $\Qbar_1$, $\Qbar_2$ the quadratic parts of
the components of $Q$. In principle these maps from $\R^4 \oplus \H$ can
contain three sorts of terms: quadratic in the first or the second
variable, or bilinear. Which terms really appear is determined by the
$J$-equivariance. Recall that $j$ acts on $\H$ by right quaternionic
multiplication and on the spaces of forms by multiplication by $-1$,
whereas $U(1)$ acts by complex multiplication on $\H$ and trivially on
the spaces of forms. This forces $\Qbar_1$ to be quadratic in the second
(quaternionic) variable and $\Qbar_2$ to be bilinear. Note that
$j$-equivariance imposes extra restrictions on these terms; clearly
$\Qbar_1 \circ j = - \Qbar_1$. The second component satisfies
$\Qbar_2 \circ j = -j \circ \Qbar_2$ if we think of $\Qbar_2$ as a linear
map $\R^4 \to \text{End}_{\C}(\H)$.

We choose a
non-degenerate $J$-invariant quadratic map $R_k \colon \H \to \R^3$ (with
the associated linear map an isomorphism, cf.~\cite{morgan-szabo:mod2})
to perturb $\Qbar_1$. For all but finitely many $\tau$, the map $\Qbar_1 +
\tau R_k$ is non-degenerate in the above sense. Admissible perturbations of
$\Qbar_2$ are of the form $(a,w) \mapsto L(a)w$, where $L(a)$ is a
$\C$-linear map which anti-commutes with the $j$-action. The space $\II$
of such maps is 4-dimensional over $\R$ and its non-zero elements are
isomorphisms. We choose the map $L_k \colon \R^4 \to \II$, $a \mapsto L_k(a)$
to be an isomorphism. Then for almost all $\tau$ the map $\Qbar_2 +\tau L_k$,
where we interpret $\Qbar_2$ as a linear map $\R^4 \to \II$, is an
isomorphism. Notice that $\Qbar_2(a,-)$ is itself an isomorphism for
$a \ne 0$; this follows from the construction of the linear perturbation
$\F$. Moreover, the norms of these linear maps are bounded from below by
$C||a||$ for some positive $C$. This means that we can choose $\tau$ small
enough so that for $a \ne  0$ the perturbation term is dominated by the
original (quadratic) map. The benefits of this perturbation are twofold;
firstly, the only solutions to the perturbed equations close to $(A_k,0,0)$
are the reducible ones. Secondly, for a generic $h \in \RR^3$, the preimage
of $(h,0,0)$ under the perturbed Kuranishi map consists of exactly one
circle of solutions, hence the point $(A_k,0,0)$ contributes $\pm 1$ to
the \sw\ invariant.

Consider now a point $(A_k,0,0)$ with $[A_k,0] \in \RR$ not $j$-fixed.
Such a point has its $j$-image in $\RR$; to make the perturbation term
$j$-equivariant in this case, we construct a $U(1)$-equivariant
perturbation at $(A_k,0,0)$ and use the $j$-action to define the
perturbation at its $j$-image. Given only $U(1)$-equivariance for
$\Qbar_1$ and $\Qbar_2$ in this case, the structure of these quadratic
maps is not so restricted. Using additional properties of the Kuranishi
map, we still conclude that $\Qbar_1$ is quadratic in the second
variable and $\Qbar_2$ is bilinear. However, the space of
$U(1)$-invariant quadratic polynomials $\H \to \R$ is four dimensional
and $\Qbar_2(a,-)$ can be any $\C$-linear map, so there is no canonical
choice of a good perturbation. To gain the same control over the solution
space as for $j$-fixed points, we endow the kernel and the cokernel with a
quaternionic structure. The perturbation terms can then be constructed
as above, using right multiplication by the quaternion $j$ in place of
the $j$-action. For the perturbation term $R_k \colon \H \to \R^3$, the
associated linear map is surjective and for all but finitely many
$\tau$, the map $\Qbar_1+ \tau R_k$ is an epimorphism in the above sense.
The perturbation of the second component gives rise to an injective map
$a \mapsto L_k(a)$; again, for all but finitely many $\tau$ the map
$a \mapsto \Qbar_2(a,-) + \tau L_k(a)$ is a monomorphism. The remark about
domination of the pertubation term $\tau L_k(a)$ by $\Qbar_2(a,-)$ holds as
above, and so do the conclusions about the solution space.

We fix a small, generic $\tau$ and define the perturbation term
as a sum of terms localized near the points of $\RR$. For a point
$[A_k,0] \in \RR$ define the perturbing map by
$$(A,\psi,w) \mapsto \tau \beta_k(A,\psi,w) \cdot (R_k(\Pi_k(\psi,w)),
L_k(P(A))\Pi_k(\psi,w),0)$$
where $\Pi_k \colon  \Gamma(W^+) \oplus \H^n  \to K_{A_k}^{\prime\prime}=\H$
is the $L^2$-orthogonal projection and $\beta_k$ is a $[0,1]$-valued function
depending smoothly on the norms of the arguments (using $A_k \equiv 0$),  that
has support inside a small neighborhood of $(A_k,0,0)$ (the  projection of
which by $(P,\Pi_k)$ is  contained in the domain of the  Kuranishi map)
and is equal to 1 on a smaller  neighborhood.  The moduli space of
solutions to the perturbed equations $\overSW=0$ still contains  the torus of
reducibles $T^*$ and it is clear from above that this torus is isolated.
\end{proof}

\section{Completion of the argument}
First we observe, following the line
of argument in \cite{morgan-szabo:mod2}, that the moduli space of solutions
to the perturbed equations $\overSW=0$ is compact.   Moreover, because the
perturbed equations we use can be connected to the unperturbed equations by a
$1$-parameter family, the count of solutions we obtain coincides with the \sw\
invariant.  In the complement of
$\Tt^*$, the action of $J$ is free, and so we can choose a small
$J$-equivariant perturbation with support away from $\Tt^*$ such that the
corresponding moduli space is smooth away from $\Tt^*$.  (This fits into
the general scheme laid down in \S4.3.6 of~\cite{donaldson-kronheimer} because
the perturbation is simply a small Fredholm section of a bundle over $(\CC'
-\Tt^*)/J$, pulled back to $\CC'$.)   Because $j$ acts freely, the
solutions in the complement of $\Tt^*$ are paired up, and this part of the
moduli space contributes an even number to the \sw\ invariant.

Along the space of reducible solutions we proceed by choosing
a small generic self-dual 2-form $\omega$ which has a nonzero harmonic
projection. If $\omega$ is  small enough, the solutions to $\overSW=
(\omega,0,0)$ in an invariant neighborhood $\UU$ of $\Tt^*$ are
described as follows. For every point in $\RR$, there is a circle of solutions
corresponding to the $U(1)$ orbit.  THere are $r$ such circles, each of which
contributes $\pm 1$ to the invariant.  All the rest of the solutions are paired
by the $j$ action, hence the statement of the theorem follows.

\section{Some homology tori}
There are a number of examples of homology tori whose \sw\
invariants one can compute directly; it is interesting to see how these are
consistent with our theorem.  The simplest are the torus $T^4$, whose \sw\
invariant is $\pm 1$, and the connected sum
$$
\#_4 S^1 \times S^3 \# \#_3 S^2 \times S^2
$$
whose \sw\ invariant vanishes.  These manifolds have determinant $1$ and $0$,
respectively.

A more interesting class of examples is the set of manifolds of the
form $X = S^1 \times M^3$, where $M$ is an orientable $3$-manifold with the
homology of a torus.  Work of Meng and Taubes shows how to compute the
invariant of $X$, in terms of the Alexander polynomial of $M$.  There are two
parts to the computation.  First, there is an identification of the \sw\
invariant of $X$ with the $3$-dimensional \sw\ invariant of $M^3$.  This is
proved by a variant of the argument proving proposition 5.1
of~\cite{morgan-szabo-taubes}.  In particular, the $\spinc$ structures on $X$
with non-vanishing \sw\ invariant all
pull back from $M$.  The main theorem of~\cite{meng-taubes} shows that the
\sw\ invariant of $M$ (and therefore of $X$) has for generating function the
multivariable Alexander polynomial of $M$.  In light of
Theorem~\ref{swtorus}, we explain how the determinant of $X$ is related to the
Alexander polynomial of $M$.

We define the determinant $\det(M)$ analogously to that of $X$, using the
$3$-fold cup
product in $H^1(M)$.  Note that the determinant of $S^1 \times M$ coincides
with that
of $M$.   The Alexander polynomial of
$M$,
$\Delta_M$, is a Laurent polynomial in variables $t_1^{\pm 1},t_2^{\pm
1},t_3^{\pm 1}$ which is defined up to multiplication by $\pm t_i$.
The relation
we need is the following:
\begin{lemma}\label{determinant} If $M$ is a homology torus, then
$$ \Delta_M(1,1,1) = \pm \det(M)^2
$$
\end{lemma}
The Lemma may be deduced from work of L.~Traldi~\cite{traldi:mubar} and
J.~Levine~\cite{levine:factorization}.  Those authors treat the Alexander
polynomial $\Delta_L$ of an $n$-component link $L$ in a homology sphere; in our
situation the homology sphere is obtained by doing surgery on a set of circles
representing a basis of
$H_1(M)$ and $L$ consists of the meridians of those circles.  If the linking
numbers between the components are all
$0$, as is the case for us, they show that
\begin{equation}\label{alex}
\frac{\Delta_L}{(t_1-1)\cdots(t_n-1)} = d_0 + \mathrm{higher\ order\ terms\ in\
} t_i-1
\end{equation}
where $d_0$ may be evaluated as a determinant involving the
$\overline\mu$-invariants of $L$.  (Compare~\cite[Corollary
1.6]{levine:factorization} and the proof of~\cite[Theorem 5.3]{traldi:mubar}.)
When there are only $3$ components, the determinant works out to be
$\overline\mu_{123}(L)^2$.  Now the quotient on the left-hand side of
equation~\eqref{alex} is the Alexander polynomial of $M$, and it has been known
for a long time~\cite{massey:linking} that the invariant
$\overline\mu_{123}(L)$
coincides with the $3$-fold Massey product.

In terms of \sw\ theory, the evaluation $\Delta_M(1,1,1)$ is the sum
of the \sw\
invariants of all of the $\spinc$ structures on $M$.  Recall that there is an
involution on the set of $\spinc$ structures, whose only fixed point is the
$\spinc$ structure $S_0$ with trivial determinant, i.e. the one we have been
studying.  Hence we have the chain of equalities and congruences
$$\SW_X(S_0) =  \SW_M(S_0) \equiv
\Delta_M(1,1,1) = \det(M)^2 \equiv \det(M) \pmod{2}
$$
which is consistent with our main theorem since $\det(M) = \det(X)$.

It is not hard to find $3$-manifolds with arbitrary determinant $\det(M)$;
a simple
construction is to take 0-framed surgery on the $n$-fold band sum of the
Borromean rings.  The case $n=2$ is illustrated below in Figure 1.  If each
copy of the Borromean rings is oriented so that the triple Massey product is
$+1$, then $\det(M) = n$.\\[2ex]
\begin{center}
\includegraphics{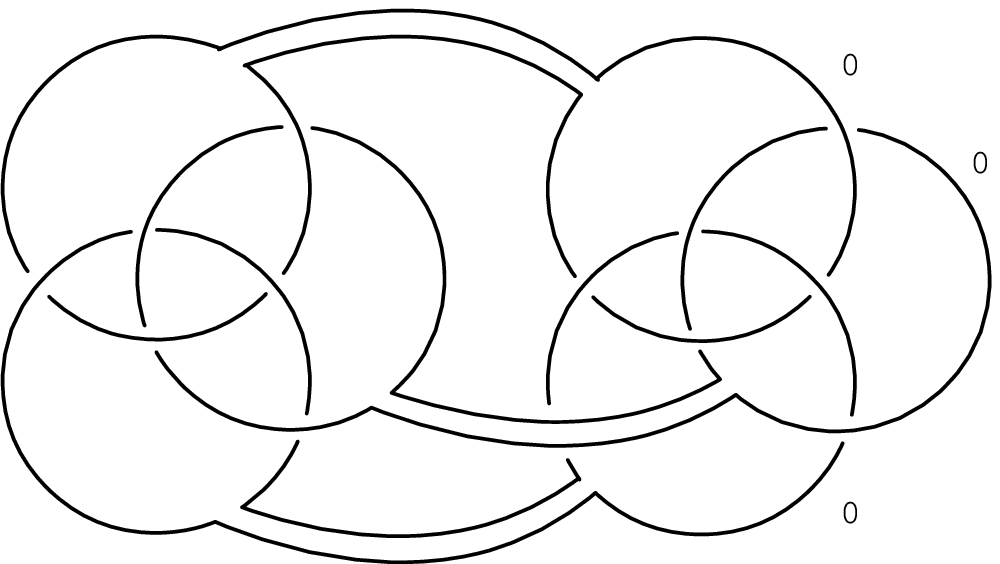}\\
Figure 1
\end{center}

\begin{remark}
The calculations assembled above give rise to a curious criterion for
a homology
torus $X^4$ to be diffeomorphic to the product of $S^1$ and a $3$-manifold.
Namely, the sum of its Seiberg-Witten invariants should be a square (up to
sign).  It would be of interest to find an example where this
criterion does not
hold, but where $X$ is homeomorphic (or perhaps homotopy equivalent) to a
product.
\end{remark}

One last class of examples is obtained via the `knot-surgery' construction of
Fintushel and Stern.  Following~\cite{fs:knots}, let $K$ be a knot in $S^3$,
with exterior $E_K$.  Remove a copy of $T^2 \times D^2$ from $T^4$, and glue in
$S^1 \times E_K$, resulting in a new manifold $X_K$ with the same cohomology as
$T^4$.  It is not hard to see that $X_K$ is in fact $S^1 \times M$, where $M$
is gotten by replacing a copy of $S^1 \times D^2 \subset T^3$ with $E_K$.  From
this, or from gluing theorems (cf.~\cite[Theorem 1.5]{fs:knots}), it follows
that the \sw\ invariant of $X_K$ is $\Delta_K(T^2)$.   To make a manifold
which is not
a product, perform this construction on three disjoint tori
$T_1,T_2,T_3$ (using
knots $K_1,K_2,K_3$) in different (non-zero) homology classes, as
in~\cite{gompf-mrowka}, to get a manifold $X_{K_1,K_2,K_3}$.
Suppose that the
knot-surgery is performed so that the circle factor in each $S^1 \times
E_{K_i}$ is glued to the same circle factor in $T^4$. The result is a product
of $S^1$ with the manifold obtained by $0$-surgery on the Borromean rings with
the knots $K_i$ tied in the three rings.    If the circle
factors in each $S^1 \times E_{K_i}$ are glued to different circles
in $T^4$ (and the knots are
non-trivial) then $X_{K_1,K_2,K_3}$ cannot be written as $S^1 $ times any
$3$-manifold.  This is verified by a fundamental group calculation; on the
other hand the \sw\ invariants are independent of the gluing and are given by
$$
\Delta_{K_1}(T_1^2)\, \Delta_{K_2}(T_2^2)\, \Delta_{K_3}(T_3^2).
$$


\end{document}